\documentclass[11pt]{article}
\newtheorem{theorem}{Theorem}
\newtheorem{lemma}[theorem]{Lemma}

 \newcommand{\nl}{\newline}
 
 \newcommand{\N}{{\bf N}}

\newcommand{\R}{{\bf R}}
 \newcommand{\C}{{\bf C}}
 
 \newcommand{\cC}{{\cal C}}
 \newcommand{\cF}{{\cal F}}

\newcommand{\dom}{{\rm Dom}}

 \newcommand{\darr}[4]{{\left\{\begin{array}{ll}
   {#1}&{#2}\\
   {#3}&{#4}
 \end{array}\right.}}

\newcommand{\sscr}[2]{\scriptstyle {#1} \atop\scriptstyle {#2}}
\newcommand{\ab}{\alpha\beta}

\title{Trace estimates and invariance of the essential spectrum}

\author{G. Barbatis}
\date{}

\begin{document}

\maketitle
\begin{abstract}
We provide sufficient conditions under which the difference of the resolvents of two higher-order
operators acting in $\R^N$ belongs to trace classes $\cC^p$. We provide explicit estimates on the norm
of the resolvent difference in terms of $L^p$ norms of the difference of the coefficients. Such inequalities are useful
in estimating the effect of localized perturbations of the coefficients.

\

\noindent {\bf AMS 2000 MSC:} 35P05 (35J30, 47F05)
\nl {\bf Keywords:}  higher-order elliptic operators, resolvent, trace classes, essential spectrum
\end{abstract}

\section{Introduction}\label{sec1}

Let $H$ and $\tilde{H}$ be second-order elliptic differential operators on $\R^N$. 
Various sufficient conditions exist under which the resolvent difference
$(\tilde{H}+1)^{-1}-(H+1)^{-1}$ is compact and, subsequently, $H$ and $\tilde{H}$
have the same essential spectrum. See \cite{H,LV,O,GG} and references therein.
These conditions typically involve some decay of the difference
$\tilde{a}_{ij}-a_{ij}$ of the respective coefficients near infinity.
Analogous results were obtained recently in \cite{GG} in a very general setting which includes
the case of higher-order operators on $\R^N$ or Laplace-Beltrami operators on manifolds.

In this note we show how an application of the Fourier
tranform can yield quantitative results of this type for higher-order self-adjoint operators of order $2m$ on $\R^N$
provided one of the operators has constant coefficients. Hence we adopt the attitude that $H$ is
is a `good', known operator and $\tilde{H}$ is a perturbed operator for which information is sought.
In our main result sufficient conditions are given under which the difference $(\tilde{H}+1)^{-1}-(H+1)^{-1}$
is not only compact on $L^2(\R^N)$ but belongs in the Schatten class $\cC^p(L^2(\R^N))$. More significantly, explicit estimates are obtained:
it is shown that if the coefficient matrix $\tilde{a}$ of $\tilde{H}$ is such
that $\tilde{a}^{-1/2}(\tilde{a}-a)a^{-1/2} \in L^p$ for some $p>N/m$ then
\begin{equation}
\|(\tilde{H}+1)^{-1}-(H+1)^{-1}\|_{\cC^p} \leq c \|\tilde{a}^{-1/2}(\tilde{a}-a)a^{-1/2}\|_{L^p}\, .
\label{radio}
\end{equation}
It is this estimate that is novel with respect to earlier work, for both the second- and higher-order case.
As a typical application, (\ref{radio}) is useful in order to estimate the effect of narrowly localised impurities
of the underlying medium; see the example following Theorem \ref{main}.
Estimates of this type were obtained by the author in \cite{B} without the assumption that one of the coefficient matrices is constant,
but the discreteness of the spectrum was a fundamental hypothesis there.

The proof uses a formula for the resolvent difference (Lemma \ref{nene}, used also in \cite{GG}) together
with a trace estimate for a class of operators acting on vector-valued fuctions. It is for the latter estimate that the Fourier
transform plays a crucial role.

We fix some notation. We work with compex-valued functions in $L^2(\R^N)$.
Given a multi-index $\alpha=(\alpha_1,\ldots,\alpha_N)$ we
use the standard notation $D^{\alpha}$ for the differential expression
$(\partial/\partial x_1)^{\alpha_1}\ldots (\partial/\partial x_N)^{\alpha_N}$.
Throughout this article we fix a positive integer $m$ and denote by $\nu=\nu(m,N)$ be the number of multi-indices
$\alpha$ of length $|\alpha|:=\alpha_1+\cdots+\alpha_N=m$. Given a vector
$v=(v_{\alpha})\in\C^{\nu}$ we denote by $v\otimes v$ the rank-one matrix
$(v_{\alpha}\overline{v}_{\beta}) \in M^{\nu\times\nu}$. The summation convention over repeated indices is adopted throughout
this article.

The $L^p$-norm of a matrix valued function $V=(V_{\ab}(x)):\R^N\to M^{\nu\times\nu}(\C)$ is defined in the standard way,
\[
\|V\|_p =\Big(\int_{\R^N}|V(x)|^pdx\Big)^{1/p},
\]
where $|V(x)|$ denotes the norm of the matrix $V(x)$ regarded as an operator on $\C^{\nu}$;
the $L^{\infty}$-norm is defined similarly.
Such a potential $V$ induces a multiplication operator on $L^2(\R^N)^{\nu}$,
also denoted by $V$, with domain $\dom(V)=\{u\in L^2(\R^N)^{\nu}\; : \; Vu\in L^2(\R^N)^{\nu}\}$, where
\[
(Vu)_{\alpha}(x) = V_{\ab}(x)u_{\beta}(x)\; , x\in\R^N \; , \quad (u_{\beta})\in L^2(\R^N)^{\nu}.
\]
If $V\in L^{\infty}$ then $V$ is a bounded operator and the two norms coincide.

We shall consider operators of order $2m$ acting on $L^2(\R^N)$ and given formally by
\label{leo}
\begin{equation}
Hu(x)=(-1)^m\sum_{\sscr{|\alpha| = m}{|\beta| = m}}
D^{\alpha}\{a_{\ab}(x)D^{\beta}u(x)\},\quad x\in\R^N\, .
\end{equation}
We assume that the complex matrix-valued function $a=(a_{\ab}(x))$ is self-adjoint and positive definite for a.e. $x\in\R^N$ and,
moreover, that
$a,a^{-1}\in L^1_{loc}(\R^N)$. To define the operator $H$ we first define the quadratic form
\[
Q(u)=\int_{\Omega}\sum_{\sscr{\mid\alpha\mid=m}{\mid\beta\mid=m}}
a_{\ab}(x)D^{\alpha}u(x)D^{\beta}\bar{u}(x)\,dx,
\]
on $C^{\infty}_c(\R^N)$ and {\em assume} that $Q$ is closable; we also denote by $Q$ its closure.
The operator $H$ is then defined as the self-adjoint operator associated with its closure. There are various sufficient conditions for the closability of $Q$, for
which we refer to \cite{Da2,GG,O,O2} and references therein.

\section{Main results}\label{sec2}

We shall initially consider uniformly elliptic operators with $a,a^{-1}\in L^{\infty}(\R^N)$ and we shall drop
this assumption in Theorem \ref{main1}.
So let $H$ be as above, with $a,a^{-1}\in L^{\infty}(\R^N)$; this implies in particular that the domain $\dom(Q)$
coincides with the Sobolev space $H^m(\R^N)$.
We define the (closed, densely defined) operator $D_m: L^2(\R^N)\to L^2(\R^N)^{\nu}$,
\[
\dom(D_m)=H^m(\R^N)\quad , \quad D_mu =(D^{\alpha}u)\, .
\]
We also denote by $b=(b_{\ab}(x))$ the square root of the matrix $a=(a_{\ab}(x))$ and define
$T=bD_m$ so that $\dom(T)=H^m(\R^N)$ and
\[H=T^*T\, .\]
We finally define the self-adjoint operator
\begin{equation}
F:  L^2(\R^N)^{\nu} \to L^2(\R^N)^{\nu}\; , \qquad F=TT^*.
\label{f}
\end{equation}
with $\dom(F)=\{ v=(v_{\alpha})\in\dom(T^*) \; : \; T^*v\in\dom(T)\}$.

Suppose now that we have two such operators $H$ and $\tilde{H}$. Keeping the above
notation and using tildes in an obvious way we have:
\begin{lemma}
There exist partial isometries $U,\tilde{U}:L^2(\R^N)\to L^2(\R^N)^{\nu}$ such that
\[
(\tilde{H}+1)^{-1}-(H+1)^{-1} =\tilde{U}^*\tilde{F}^{1/2}(\tilde{F}+1)^{-1}\tilde{a}^{-1/2}(a-\tilde{a})
a^{-1/2}(F+1)^{-1}F^{1/2}U.
\]
\label{nene}
\end{lemma}
{\em Proof.}
We write the identity \cite[p271]{De}
\[
(S^*S+1)^{-1}+S^*(SS^*+1)^{-1}S=I
\]
first for $S=T$, then for $S=\tilde{T}$ and subtract the two relations; we obtain
\begin{eqnarray*}
(\tilde{H}+1)^{-1}-(H+1)^{-1}&=&-D_m^*\tilde{a}^{1/2}(\tilde{F}+1)^{-1}\tilde{a}^{1/2}D_m
+D_m^* a^{1/2}(F+1)^{-1}a^{1/2}D_m \\
&=&-D_m^*\Big\{  (D_mD_m^*+\tilde{a}^{-1})^{-1} -(D_mD_m^*+a^{-1})^{-1}\Big\}D_m \\
&=&D_m^*(D_mD_m^*+\tilde{a}^{-1})^{-1}(\tilde{a}^{-1}-a^{-1})(D_mD_m^*+a^{-1})^{-1}D_m\\
&=&D_m^*\tilde{a}^{1/2}(\tilde{F}+1)^{-1}\tilde{a}^{1/2}(\tilde{a}^{-1}-a^{-1})a^{1/2}(F+1)^{-1}a^{1/2}D_m\\
&=&D_m^*\tilde{a}^{1/2}(\tilde{F}+1)^{-1}\tilde{a}^{-1/2}(a-\tilde{a})a^{-1/2}(F+1)^{-1}a^{1/2}D_m\, .
\end{eqnarray*}
Using polar decomposition, we write
\[
a^{1/2}D_m=F^{1/2}U \;\; , \qquad  \tilde{a}^{1/2}D_m=\tilde{F}^{1/2}\tilde{U} \, ,
\]
where $U,\tilde{U}: L^2(\R^N)\to L^2(\R^N)^{\nu}$ are partial isometries; this completes the proof.
$\hfill //$

{\bf Remark.} It is an immediate consequence of Lemma \ref{nene} that
\[
\|(\tilde{H}+1)^{-1}-(H+1)^{-1}\| \leq \frac{1}{4}\|\tilde{a}^{-1/2}(\tilde{a}-a)a^{-1/2}\|_{\infty}\, .
\]
This of course does not contain any information on possible compactness of the resolvent difference;
nor is it any useful in the context of the example following Theorem \ref{main}. In what follows we shall
concentrate with the case where $\tilde{a}^{-1/2}(\tilde{a}-a)a^{-1/2}\in L^p(\R^N)$ for some $p<\infty$.

To proceed we define the weighted $L^p$ spaces
\[
L^p( \R_+ ,t^{\frac{N-2m}{2m}}dt)= \Big\{ g: \R_+\to\R \; : \;  \int_0^{\infty}|g(t)|^pt^{\frac{N-2m}{2m}}dt <+\infty \Big\}\; ,
\quad 1\leq p<\infty\, ,
\]
equipped with the natural norm, which for simplicity we denote by $\|\cdot\|_p^{*}\,$.
The next lemma is a vector-valued version of \cite[Theorem 4.1]{S}; we note that it does not follow directly
from that result. Although the proof is very similar to that of \cite{S}, we present it for the sake of
complete
\begin{lemma}
Assume that the operator $H$ has constant coefficients.
Let $V=(V_{\ab}(x))$ be a matrix-valued function and let $g:[0,+\infty)\to\R$
be a bounded continuous function with $g(0)=0$. If for some
$p\in [2,+\infty)$ there holds $V\in L^p(\R^N)$ and $g\in L^p( \R_+ ,t^{\frac{N-2m}{2m}}dt)$, then
$Vg(F)\in \cC^p(L^2(\R^N)^{\nu})$ and moreover
\begin{equation}
\|Vg(F)\|_{\cC^p} \leq c^{1/p}\|V\|_{L^p}\|g\|_p^*\, .
\label{sis}
\end{equation}
for a constant $c=c(H)$.
\label{leve}
\end{lemma}
{\em Proof.} We may assume that both $V$ and $g$ have compact supports since the general case will then follow by approximation
and an application of the dominated convergence theorem for trace ideals. Also, it is enough to establish (\ref{sis}) for $p=2$
and $p=\infty$ since the intermediate cases will then follow by interpolation \cite[Theorem 2.9]{S}.

Let us denote by $\cF$ the Fourier transform regarded as a unitary operator on $L^2(\R^N)$; we use the same symbol for the unitary operator induced
component-wise on $L^2(\R^N)^{\nu}$. The fact that the differential operator $F$ (cf. (\ref{f})) has constant coefficients implies that $F$ is unitarily
equivalent via $\cF$ to a multiplication operator on $L^2(\R^N)^{\nu}$.
More precisely, for $\xi\in\R^N$ (the variable in the Fourier space) let us define the vector
\[  
(B(\xi))_{\alpha}=b_{\alpha\gamma}\xi^{\gamma}.
\]
We then have $\cF F\cF^{-1}v=B(\xi)\otimes B(\xi)$, that is
\[ (\cF F\cF^{-1}v)_{\alpha}(\xi)= (B(\xi)\otimes B(\xi))_{\ab}v_{\beta}(\xi) =
b_{\alpha\gamma}\overline{b}_{\beta\delta}\xi^{\gamma+\delta}v_{\beta}(\xi).
\]
This implies that
\[
\cF F^k\cF^{-1}= |B(\xi)|^{2k-2}B(\xi)\otimes B(\xi)\; , \quad k=1,2,\ldots\, ,
\]
and hence
\[
\cF p(F)\cF^{-1}=p(0)I  + \Big[p(|B(\xi)|^2) -p(0)I\Big] |B(\xi)|^{-2}B(\xi)\otimes B(\xi)
\]
for any polynomial $p(\cdot)$. Direct computation then shows that if $p(\cdot)$ does not vanish on $[0,\infty)$ then
\[
\cF p^{-1}(F)\cF^{-1}=p^{-1}(0)I  + \Big[p^{-1}(|B(\xi)|^2) -p^{-1}(0)I\Big] |B(\xi)|^{-2}B(\xi)\otimes B(\xi).
\]
Compactifying $[0,\infty)$ we obtain by an application of the Stone-Weierstrass theorem that 
\[ \cF g(F) \cF^{-1} =g(|B(\xi)|^2)|B(\xi)|^{-2}B(\xi)\otimes B(\xi)=: (Lg)(\xi)\, , \quad \xi\in\R^N, \]
for any continuous function $g$ on $[0,+\infty)$ with $g(0)=0$. We note that $L$ is a linear map from the space of all
such $g$ to the space of matrix-valued functions and $Lg$ is a multiplication operator in $L^2(\R^N)^{\nu}$.

It follows that $g(F)$ has an integral kernel depending only on $x-y$, $k_g=k_g(x-y)$, where
\[
k_g(x) =(2\pi)^{-N}(\cF^{-1}Lg)(x)\, , \qquad x\in\R^N\, ;
\]
hence $Vg(F)$ has the integral kernel $V(x)k_g(x-y)$.
It follows that $Vg(F)$ is a Hilbert-Schmidt operator with Hilbert-Schmidt norm given by
\begin{eqnarray}
\|Vg(F)\|_{\cC^2}^2&=&\int_{R^N}\int_{R^N}|V(x)k_g(x-y)|^2dx\, dy \nonumber \\
&\leq&\int_{R^N}\int_{R^N}|V(x)|^2  |k_g(x-y)|^2dx\, dy \nonumber \\
&=&\|V\|_{L^2}^2\|k_g\|_2^2 \label{kik}\\
&=&(2\pi)^{-2N}\|V\|_{L^2}^2\|Lg\|_2^2 \, .
\nonumber
\end{eqnarray}
Using the homogeneity of the symbol $A(\xi)$ of $H$ we obtain by an application of the coarea
formula,
\begin{eqnarray}
\|Lg\|_2^2 &=&\int_{\R^N}g^2(|B(\xi)|^2)|B(\xi)|^{-4}|B(\xi)\otimes B(\xi)|^2 d\xi \nonumber \\
&=&\int_{\R^N}g^2(|B(\xi)|^2)d\xi  \label{dim}\\
&\leq&c(H)\int_0^{\infty}g^2(t)t^{\frac{N-2m}{2m}}dt  \nonumber\\
&=&c(H)\|g\|_2^{*2}\, . \nonumber
\end{eqnarray}
Combining (\ref{kik}) and (\ref{dim}) we conclude that
\[
\| Vg(F)\|_{\cC^2}^2 \leq (2\pi)^{-2N}c(H)\|V\|_{L^2}^2\|g\|_2^{*2} \, .
\]
We also have
\[
\| Vg(F)\|_{\cC^{\infty}}\leq \|V\|_{L^{\infty}}\|g\|_{\infty};
\]
here we note that the compactness of $Vg(F)$ follows from our assumption on the supports of $V$, $g$, which,
by our argument above, implies that $V$ and $g$ are both $L^2$ and hence that $Vg(F)$ is Hilbert-Schmidt. 
This completes the proof of the lemma. $\hfill //$

\begin{theorem}
Let $H$ and $\tilde{H}$ be uniformly elliptic self-adjoint operators of order $2m$ and let $a$ and $\tilde{a}$ be
the respective coefficient matrices. Assume that $H$ has constant coefficients. Then for any $p\in (N/m ,\infty)$
there exists a positive constant $c=c(p,H)$ such that
\[
\|(\tilde{H}+1)^{-1}-(H+1)^{-1}\|_{\cC^p} \leq c\|\tilde{a}^{-1/2}(\tilde{a}-a)a^{-1/2}\|_{L^p} \, .
\]
\label{main}
\end{theorem}
{\em Proof.} Setting $g(t)=t^{1/2}(t+1)^{-1}$ and using Lemmas \ref{nene} and \ref{leve} we obtain
\begin{eqnarray*}
\|(\tilde{H}+1)^{-1}-(H+1)^{-1}\|_{\cC^p}&=&\|g(\tilde{F})\tilde{a}^{-1/2}(\tilde{a}-a)
a^{-1/2}g(F)\|_{\cC^p}\\
&\leq&\frac{1}{2}\|\tilde{a}^{-1/2}(\tilde{a}-a)a^{-1/2}g(F)\|_{\cC^p} \\
&\leq&\frac{1}{2}c^{1/p}\|\tilde{a}^{-1/2}(\tilde{a}-a)a^{-1/2}\|_{L^p}\|g\|_p^*.
\end{eqnarray*}
The proof is completed if we note that $\|g\|_p^*<\infty$ if and only if $p>N/m$. $\hfill //$

{\bf Example.} Suppose that $H$ is an operator with constant coefficients $a=\{a_{\ab}\}$ as above describing some physical phenomenon
and assume that the presence of some localized impurities on a set $U$ of finite volume yields a new coefficient matrix,
\[
\tilde{a}(x)=\darr{a+b(x),}{x\in U,}{a,}{x\not\in U\, ,}
\]
where $b\in L^{\infty}(U)$. We then have
\[ \|\tilde{a}-a\|_p =\|b\|_{L^p(U)}\leq \|b\|_{L^{\infty}(U)}|U|^{1/p}.\]
Hence we have a precise esimate on the effect of the given impurity in terms of the volume $|U|$.

Finally, at the cost of having at the left-hand side the operator norm instead of a $\cC^p$ norm, we drop the uniform ellipticity assumption on $\tilde{H}$.
\begin{theorem}
Let $H$ and $\tilde{H}$ be self-adjoint elliptic operators of order $2m$ and let $a$ and $\tilde{a}$ be
the respective coefficient matrices. Assume that $H$ has constant coefficients. If 
$\tilde{a}^{-1/2}(\tilde{a}-a) \in L^p$ for some $p\in (N/m,\infty)$ then there exists a positive constant $c=c(p,H)$ such that
\begin{equation}
\|(\tilde{H}+1)^{-1}-(H+1)^{-1}\|  \leq c\|\tilde{a}^{-1/2}(\tilde{a}-a)a^{-1/2}\|_{L^p} \, .
\label{chess}
\end{equation}
\label{main1}
\end{theorem}
{\em Proof.} Using the diagonalization of $\tilde{a}(x)$ we define for each $x\in\R^N$ the matrices
\[
\tilde{a}_n(x) =\max\Big\{ 1/n , \min\{ \tilde{a}(x) , n\}\Big\}\, , \quad n=1,2,\ldots \, .
\]
The corresponding operators $\tilde{H}_n$ are then uniformly elliptic and Theorem \ref{main} implies that for $p>N/m$
there exists $c=c(p,H)$ such that
\begin{equation}
\|(\tilde{H}_n+1)^{-1}-(H+1)^{-1}\| \leq c\|\tilde{a}_n^{-1/2}(\tilde{a}_n-a)a^{1/2}\|_{L^p} \, ,
\label{q1}
\end{equation}
where the constant $c$ is independent of $n\in\N$. Now, it follows from \cite[p118]{Da1} and \cite[Theorem 1.2.3]{Da2}
that $(\tilde{H}_n+1)^{-1}\to (\tilde{H}+1)^{-1}$ strongly as $n\to +\infty$. This together with
(\ref{q1}) and Lebesgue's dominated convergence theorem yields (\ref{chess}). $\hfill //$

{\bf Remark.} A version of Lemma \ref{leve} for operators with variable coefficients would extend our results to
the case where both $H$ and $\tilde{H}$ have variable coefficients. This is an open problem.

{\bf Acknowledgment.} The material for this article has been essentially extracted from the author's Ph.D.
thesis. Hence I thank E.B. Davies once again for his help and guidance when this work was being carried out.
I also thank the referee for useful comments.



G. Barbatis \\  
Department of Mathematics \\
University of Ioannina\\ 
45110 Ioannina \\
Greece

\end{document}